\documentclass[matann2]{svjour}

\usepackage{amsmath}
\usepackage{amsfonts}
\usepackage{latexsym}
\usepackage{mathrsfs}
\usepackage{times}
\usepackage{array}
\usepackage{amscd}
\usepackage{a4}
\usepackage[mathcal]{euscript}
\usepackage{amssymb}
\usepackage{pstricks}

\input{boxedeps} 
\SetEPSFDirectory{} 
\HideDisplacementBoxes
\SetRokickiEPSFSpecial  
\DeclareMathAlphabet{\ams}{U}{msb}{m}{n}

\def\sl{\text{SL}}

\def\gal{\text{Gal}}

\def\rk{\text{rk}\,}

\def\HH{\mathcal H}

\def\LL{\mathcal L}

\def\ve{\varepsilon}

\def\aa{\alpha}

\def\vphi{\varphi}

\def\ve{\varepsilon}

\def\Z{\ams{Z}}
\def\R{\ams{R}}

\def\quo{/\kern -.45em\sim}

\newpsobject{showgrid}{psgrid}{subgriddiv=1,griddots=5,gridlabels=6pt}
\title{Graphs, free groups and the Hanna Neumann conjecture}

\author{Brent Everitt\thanks{Some of the results of this paper were 
obtained while the author was visiting
the Department of Mathematics, University of Adelaide, Australia. 
He is grateful for their hospitality. 
He is also grateful to the referee, whose careful reading and
numerous suggestions have significantly improved (and shortened) the exposition
of this paper.
}
}

\institute{\textsc{Brent Everitt}:
Department of Mathematics, University of York, York
YO10 5DD, United Kingdom. \email{bje1@york.ac.uk}
}

\titlerunning{}
\authorrunning{Brent Everitt 
}

\begin{document}

\maketitle



\begin{abstract}
A new bound for the rank of the intersection of finitely generated subgroups
of a free group is given, formulated in topological terms, and
very much in the spirit of Stallings \cite{Stallings83}. The bound is a 
contribution to 
the strengthened
Hanna Neumann conjecture.
\end{abstract}


\section*{Introduction}

This paper is about the interplay between graphs and free groups, with
particular application to subgroups of free groups. This subject
has a long history, where one approach is to treat graphs as purely combinatorial
objects (as in for instance \cite{Imrich76,Tardos96,Tardos92}), while 
another (for example \cite{Scott79}), is to treat them topologically
by working in the category of $1$-dimensional CW complexes.

We prefer a middle way, where 
to quote Stallings \cite{Stallings83} who initiated it,
graphs are ``something purely combinatorial
or algebraic'', but also
one may apply to them topological machinery, motivated by their geometrical
realizations. We use this
to give a new bound for the rank of the intersection of two
finitely generated subgroups of a free group (Theorems 
\ref{pullback:rankestimate} and 
\ref{algebraic:shn}), and to formulate graph theoretic versions of some other
classical results. The first section sets up the combinatorial-topological background;
\S\ref{section:invariants} studies graphs of finite rank; the topological
meat of the paper is \S\ref{section:pullbacks} and the group theoretic
consequences explored in \S\ref{free}.

\section{Preliminaries from the topology of graphs}\label{section:topological}

%

A {\em combinatorial $1$-complex\/} or {\em graph\/} \cite[\S 1.1]{Gersten83} 
(see also \cite{Cohen89,Collins98,Serre03,Stallings83})
is 
a set $\Gamma$ with involutary $^{-1}:\Gamma\rightarrow\Gamma$
and idempotent $s:\Gamma\rightarrow V_\Gamma$, (ie: $s^2=s$) maps,
where $V_\Gamma$ is the
set of fixed points of $^{-1}$. Thus a graph has {\em vertices\/} $V_\Gamma$
and {\em edges\/} $E_\Gamma:=\Gamma\setminus V_\Gamma$ with
(i). $s(v)=v$ for all $v\in V_\Gamma$;
(ii). $v^{-1}=v$ for all $v\in V_\Gamma$, $e^{-1}\in E_\Gamma$ and 
$e^{-1}\not= e=(e^{-1})^{-1}$ for all 
$e\in E_\Gamma$. 
The edge $e$ has start vertex $s(e)$ and terminal vertex
$t(e):=s(e^{-1})$; an arc is an edge/inverse
edge pair;
a pointed graph is a pair
$\Gamma_v:=(\Gamma,v)$ for $v\in\Gamma$ a vertex.

A map of graphs is a set map $f :\Gamma\rightarrow\Lambda$
with $f(V_\Gamma)\subseteq V_\Lambda$
that commutes with $s$ and $^{-1}$, and preserves
dimension if $f(E_\Gamma)\subseteq E_\Lambda$. An isomorphism
is a dimension preserving map, bijective on the vertices and edges.
A map $f:\Gamma_v\rightarrow\Lambda_u$ of pointed graphs is a graph map
$f:\Gamma\rightarrow\Lambda$ with $f(v)=u$.

A graph $\Gamma$
has a functorial geometric realization
as a $1$-dimensional CW complex $B\Gamma$
(see, eg: \cite[\S 1.3]{Gersten83}) with a graph map $f:\Gamma\rightarrow\Lambda$
inducing a regular cellular map $Bf:B\Gamma\rightarrow B\Lambda$ 
of CW complexes, in the sense of \cite[\S4]{Lundell69}. 
Thus, one may transfer to graphs and their maps topological notions
and adjectives
(connected, fundamental group, homology, covering map, etc...)
from their geometrical realizations.

If $\Lambda\hookrightarrow\Gamma$ is a subgraph,
we will write $\Gamma/\Lambda$ for the resulting quotient graph
and quotient map $q:\Gamma\rightarrow\Gamma/\Lambda$.
For a set
$\Lambda_i\hookrightarrow\Gamma$, $(i\in I)$,
of mutually disjoint subgraphs, we will write $\Gamma/\Lambda_i$ for the graph
resulting from taking successive quotients by the $\Lambda_i$. 
The coboundary $\delta\Lambda$ of a subgraph consists of those 
edges $e\in\Gamma$ with $s(e)\in\Lambda$ and $t(e)\not\in\Lambda$;
equivalently, it is those 
edges $e\in\Gamma$ with
$sq(e)$ the vertex $q(\Lambda)$
in the quotient graph
$q:\Gamma\rightarrow\Gamma/\Lambda$.
The real line graph $\R$ has vertices $V_\R=\{v_k\}_{k\in\Z}$ and
edges $E_\R=\{e_k^{\pm 1}\}_{k\in\Z}$ with $s(e_k)=v_k,s(e_k^{-1})=v_{k+1}$.

We have the obvious notion of path and in particular, a spur
is a path that successively traverses both edges of an arc, and a path
is reduced when it contains no spurs. 
A tree is a simply connected graph and a forest a graph, all of whose
connected components are trees. Any connected graph has a spanning tree
$T\hookrightarrow\Gamma$ with the homology $H_1(\Gamma)$ free abelian 
on the set of \emph{arcs\/} of $\Gamma$ omitted by $T$, and the
rank $\rk\Gamma$ of $\Gamma$ (connected) defined to be $\rk_\Z H_1(\Gamma)$.
If $\Gamma$ has finite rank then $\rk\Gamma-1=-\chi(\Gamma)$, and
if $\Gamma$ is finite, locally finite, connected, then
$2(\rk\Gamma-1)=
|E_\Gamma|-2|V_\Gamma|$.
If $\Gamma$ is connected and $T_i\hookrightarrow\Gamma$ a set of 
mutually disjoint trees, then 
the fundamental group is unaffected by their excision:
$\pi_1(\Gamma,v)\cong\pi_1(\Gamma/T_i,q(v))$ and so
$\rk\Gamma=\rk\Gamma/T_i$.

If $\Lambda$ is a connected graph and $v$ a vertex, 
then the spine $\widehat{\Lambda}_v$ of $\Lambda$ at $v$, 
is defined
to be the union in $\Lambda$ of all closed reduced paths starting at $v$.
It is easy to show that
$\widehat{\Lambda}_v$ is connected with 
$\rk\widehat{\Lambda}_v=\rk\Lambda$, that
every closed reduced path
starting at $u\in\widehat{\Lambda}_v$ is contained in $\widehat{\Lambda}_v$,
and
an isomorphism $\Lambda_u\rightarrow\Delta_v$
restricts to a isomorphism $\widehat{\Lambda}_u\rightarrow\widehat{\Delta}_v$
(so that spines are invariants of graphs).

\parshape=8
0pt\hsize 0pt\hsize 
0pt.76\hsize 0pt.76 \hsize 0pt.76\hsize 0pt.76\hsize 0pt.76\hsize 
0pt\hsize 
If $\Lambda_1,\Lambda_2$ and $\Delta$ are graphs and 
$f_i:\Lambda_i\rightarrow \Delta$ 
maps of graphs, then
the pullback $\Lambda_1\prod_\Delta \Lambda_2$  has vertices
(resp. edges) the $x_1\times x_2$, $x_i\in V_{\Lambda_i}$ (resp.
$x_i\in E_{\Lambda_i}$) such that $f_1(x_1)=f_2(x_2)$, and
$s(x_1\times x_2)=s(x_1)\times s(x_2)$, $(x_1\times x_2)^{-1}
=x_1^{-1}\times x_2^{-1}$ (see \cite[page 552]{Stallings83}). 
Taking $\Delta$ to be the trivial graph 
gives the product
$\Lambda_1\prod \Lambda_2$. 
Define maps $t_i:\Lambda_1\prod_\Delta\Lambda_2\rightarrow\Lambda_i$ to be
the compositions $\Lambda_1\prod_\Delta\Lambda_2\hookrightarrow
\Lambda_1\prod\Lambda_2\rightarrow\Lambda_i$, with the second map the projection
$x_1\times x_2\mapsto x_i$. Then the
$t_i$ are 
dimension preserving maps
making the diagram
commute, and the pullback is universal with this property.
\vadjust{\hfill\smash{\lower 1mm
\llap{
\begin{pspicture}(2,2)
\rput(1.3,-1.15){
\rput(-1.5,1.5){
\rput(0,2){$\Lambda_1\prod_\Delta \Lambda_2$}
\rput(0,0){$\Lambda_1$}
\rput(2,2){$\Lambda_2$}
\rput(2,0){$\Delta$}
\psline[linewidth=.1mm]{->}(0,1.7)(0,.3)
\psline[linewidth=.1mm]{->}(.9,2)(1.7,2)
\psline[linewidth=.1mm]{->}(.3,0)(1.7,0)
\psline[linewidth=.1mm]{->}(2,1.7)(2,.3)
\rput(.25,1){$t_1$}
\rput(1.3,1.8){$t_2$}
\rput(1,.25){$f_1$}
\rput(1.75,1){$f_2$}
}
}
\end{pspicture}
}}}\ignorespaces

In general the pullback need not be connected, but if
the $f_i:\Lambda_{u_i}\rightarrow\Delta_v$ are pointed maps then 
the \emph{pointed\/} pullback $(\Lambda_1\prod_\Delta\Lambda_2)_{u_1\times u_2}$
is the connected component of the pullback containing the vertex $u_1\times u_2$
(and we then have a pointed version of 
the diagram above).

There is a ``co''-construction, the pushout, for
dimension preserving maps of graphs,
$f_i:\Delta\rightarrow \Lambda_i$,
although it will play a lesser role for us (see \cite[page 552]{Stallings83}).
The principal example for us is
the wedge sum $\Lambda_1\bigvee_{\Delta}\Lambda_2$.

Graph coverings $f:\Lambda\rightarrow\Delta$
can be characterized combinatorially as dimension preserving maps
such that
for every vertex $v\in \Lambda$, $f$ is a bijection from the set of 
edges in $\Lambda$
with start vertex $v$ to the set of edges in $\Delta$ with start
vertex $f(v)$.
Graph coverings have the usual path and homotopy lifting properties
\cite[\S 4]{Stallings83}, and
from now on, all coverings will be maps between connected complexes
unless stated otherwise, and we will write $\text{deg}(\Lambda\rightarrow\Delta)$
for the degree of the covering.
A covering is Galois
if for all closed paths $\gamma$ at
$v$, the lifts of $\gamma$ to each vertex of 
the fiber of $v$
are either all closed or all non-closed.

\begin{proposition}\label{topology:coverings:result400}
Let $\Lambda$ be a graph and $\Upsilon_1,\Upsilon_2\hookrightarrow\Lambda$
subgraphs of the form,
$$
\begin{pspicture}(14,1)
\rput(4,0){
\rput(3,.5){\BoxedEPSF{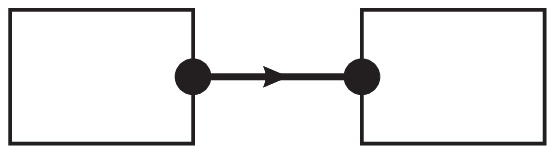 scaled 750}}
\rput(.5,.5){$\Lambda=$}
\rput(1.6,.5){$\Upsilon_1$}\rput(4.4,.5){$\Upsilon_2$}
\rput(3,.7){$e$}
}
\end{pspicture}
$$
\noindent(i). If $f:\Lambda\rightarrow\Delta$ is a covering with $\Delta$ 
single vertexed, then the real line is a subgraph 
$g:\R\hookrightarrow\Lambda$,
with $g(e_0)=e$ and $fg(e_k)=f(e)$ for all $k\in\Z$.
%
(ii). If $\Upsilon_1$ is a tree, $\Lambda\rightarrow\Delta$ and
$\Gamma\rightarrow\Delta$ coverings, and
$\Upsilon_2\hookrightarrow\Gamma$ a subgraph,
then
there is an intermediate
covering $\Lambda{\rightarrow}\Gamma{\rightarrow}\Delta$.
%
(iii). If $\Upsilon_1$ is a tree, and $\Psi\rightarrow\Lambda$ a covering,
then $\Psi$ has the same form as $\Lambda$ for some subgraphs 
$\Upsilon'_1,\Upsilon'_2\hookrightarrow\Psi$ and with $\Upsilon'_1$ a tree.
\end{proposition}

\begin{proof}
These are easy exercises using path lifting. For (i), 
build $\R\hookrightarrow\Lambda$ by taking successive lifts of the
edge $f(e)\in\Delta$.
For (ii), it suffices
to find a map $\Lambda\rightarrow\Gamma$
commuting with the two coverings given. Let it coincide with 
$\Upsilon_2\hookrightarrow\Gamma$ on $\Upsilon_2$, and on 
$\Upsilon_1$, project to $\Delta$ and then lift to $\Gamma$.
For (iii), take $\Upsilon_1'$ to be the union of lifts of reduced paths
from $t(e)$ to the vertices of $\Upsilon_1$.
\qed
\end{proof}

If $f:\Lambda\rightarrow\Delta$ is a covering 
and $T\hookrightarrow\Delta$ a 
tree, then
path and homotopy lifting give that
$f^{-1}(T)$ is a forest such that if 
$T_i\hookrightarrow\Lambda, (i\in I)$ are the 
component trees, then $f$ maps each $T_i$
isomorphically onto $T$. There is then an induced
covering $f':\Lambda/T_i\rightarrow\Delta/T$, 
defined by $f'q'=qf$
where $q,q'$ are the quotient maps, and such that
$\text{deg}(\Lambda/T_i\rightarrow\Delta/T)=
\text{deg}(\Lambda\rightarrow\Delta)$.

If $f:\Lambda_u\rightarrow\Delta_v$ is a covering
then intermediate coverings 
$\Lambda_u{\rightarrow}\Gamma_x{\rightarrow}\Delta_v$
and
$\Lambda_u{\rightarrow}\Upsilon_{y}{\rightarrow}\Delta_v$
are equivalent 
if and only if there is a isomorphism $\Gamma_x\rightarrow\Upsilon_{y}$
making the obvious diagram commute. 
Then the set
$\mathcal{L}(\Lambda_u,\Delta_v)$ of equivalence classes
of intermediate coverings
is a lattice 
with join $\Gamma_{x_1}\vee\Upsilon_{x_2}$ the pullback 
$(\Gamma\prod_\Delta\Upsilon)_{x_1\times x_2}$, 
meet $\Gamma_{x_1}\wedge\Upsilon_{x_2}$
the pushout $(\Gamma\coprod_\Lambda\Upsilon)_{g(x_i)}$ 
($g$ the covering $\Lambda_u\rightarrow\Gamma_{x_1}$), 
a $\widehat{0}=\Delta_v$
and a $\widehat{1}=\Lambda_u$.
The incessant pointing of  covers is annoying, but essential if one 
wishes to work with
{\em connected\/} intermediate coverings and also have a lattice structure 
(both of which
we do). The problem is the pullback: because it is not in general connected, 
we need the
pointing to tell us which component to choose.



We end 
the preliminaries by observing 
that the excision of trees has 
little effect
on the lattice $\mathcal{L}(\Lambda,\Delta)$.
Let $f:\Lambda_u\rightarrow\Delta_v$ be a covering, $T\hookrightarrow\Delta$ a 
spanning tree, $T_i\hookrightarrow\Lambda$ the components of $f^{-1}(T)$, and
$f:(\Lambda/T_i)_{q(u)}\rightarrow(\Delta/T)_{q(v)}$ the induced covering (where
we have (ab)used $q$ for both quotients and $f$ for both coverings).
One can then show (either by brute force, or using the 
Galois correspondence between $\mathcal{L}(\Lambda,\Delta)$ and the
subgroup lattice of the group $\gal(\Lambda,\Delta)$
of covering transformations), that
there is a degree and rank preserving isomorphism of lattices 
$\LL(\Lambda,
\Delta)\rightarrow
\LL(\Lambda/T_i,
\Delta/T)$,
that sends Galois coverings to Galois coverings, and the equivalence class of 
$\Lambda_u\rightarrow\Gamma_{x}\stackrel{r}{\rightarrow}\Delta_v$ to the 
equivalence class of 
$\Lambda/T_i\rightarrow\Gamma/T'_i\rightarrow\Delta/T$,
with $T'_i\hookrightarrow\Gamma$ the components of $r^{-1}(T)$.
We will call this process \emph{lattice excision\/}.

\section{Graphs of finite rank}
\label{section:invariants}

This section is devoted to a more detailed study of the
coverings $\Lambda\rightarrow\Delta$ where $\rk\Lambda<\infty$.

\begin{proposition}
\label{finite_rank_characterisation200}
Let $\Lambda$ be a connected graph, $\Gamma\hookrightarrow\Lambda$ a 
connected subgraph and $v\in\Gamma$ a vertex such that every closed reduced
path at $v$ in $\Lambda$ is contained in $\Gamma$. Then $\Lambda$ has a 
wedge sum decomposition $\Lambda=\Gamma\bigvee_\Theta\Phi$
with $\Phi$ a forest 
and no two vertices of the image of 
$\Theta\hookrightarrow\Phi$ lying in the same component.
\end{proposition}

\begin{proof}
Consider an edge $e$ of $\Lambda\setminus\Gamma$ having 
at least 
one of its end vertices $s(e)$ or
$t(e)$, in $\Gamma$. 
For definiteness we can assume, by relabeling
the edges in the arc containing $e$, that it is $s(e)$ that is a vertex of 
$\Gamma$. If $t(e)\in\Gamma$ 
then by traversing a reduced path in $\Gamma$ from $v$ to $s(e)$, crossing
$e$ and a reduced path in $\Gamma$ from $t(e)$ to $v$, we get a closed 
reduced path not contained in $\Gamma$, a contradiction. 
Thus $t(e)\not\in\Gamma$. Let $T_e$ be the union 
of all the reduced paths
in $\Lambda\setminus\{e\}$ starting at $t(e)$, so we have the situation
as in (a):
$$
\begin{pspicture}(0,0)(13,3)
\rput(1.5,1.7){\BoxedEPSF{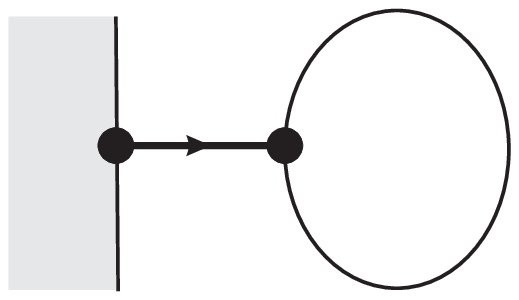 scaled 650}}
\rput(0.1,1.7){$\Gamma$}
\rput(2.35,1.7){$T_e$}\rput(1.1,1.95){$e$}
\rput(1.5,0){(a)}
\rput(6.5,1.7){\BoxedEPSF{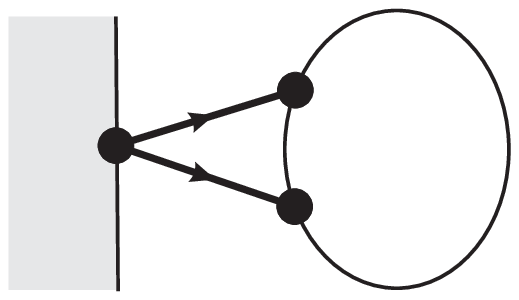 scaled 650}}
\rput(5.1,1.7){$\Gamma$}
\rput(7.4,1.7){$T_e$}
\rput(6.1,2.1){$e$}\rput(6.1,1.3){$e'$}
\rput(6.5,0){(b)}
\rput(11.5,1.7){\BoxedEPSF{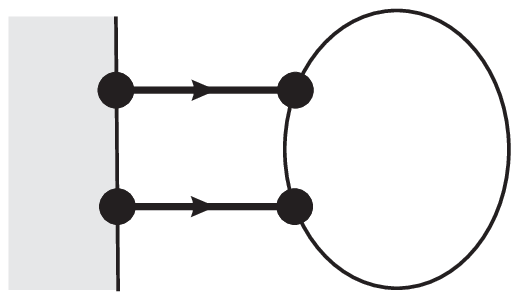 scaled 650}}
\rput(10.1,1.7){$\Gamma$}
\rput(12.4,1.7){$T_e$}
\rput(11.1,2.3){$e$}\rput(11.1,1.1){$e'$}
\rput(11.5,0){(c)}
\end{pspicture}
$$
If $\gamma$ is a non-trivial closed path in
$T_e$ starting at $t(e)$, then a path from $v$ to $t(e)$, 
traversing $\gamma$, and going the same way back to
$v$ cannot be reduced. But the only place a spur can occur is in $\gamma$
and so $T_e$ is a tree.
If $e'$ is another edge of $\Lambda\setminus\Gamma$ with
$s(e')\in\Gamma$ then we claim that neither of the 
two situations (b) and (c) above can occur, ie: $t(e')$ is not a vertex of 
$T_e$. For otherwise, a reduced closed path in $T_e$ from $t(e)$ to
$t(e')$ will give a reduced closed path at $v$ not in $\Gamma$.
Thus, another edge $e'$ yields a tree $T_{e'}$ defined like
$T_e$, but disjoint from it. Each component of $\Phi$ is thus obtained 
this way.
\qed
\end{proof}

\begin{corollary}\label{invariants:result400}
$\Lambda$ connected is of finite rank if and only if for any vertex $v$, the
spine $\widehat{\Lambda}_v$ is finite, locally finite.
\end{corollary}

\begin{proof}
Proposition \ref{finite_rank_characterisation200}
gives the wedge sum decomposition 
$\Lambda=\widehat{\Lambda}_v\bigvee_\Theta\Phi$, and by connectedness, any spanning tree
$T\hookrightarrow\Lambda$ must contain the forest $\Phi$ as a subgraph. Thus if $\Lambda$
has finite rank, then $\widehat{\Lambda}_v$ is a tree with finitely many edges added,
hence finite.
Conversely, a finite spine has
finite rank and $\rk\Lambda=\rk\widehat{\Lambda}_v$.
\qed
\end{proof}

Thus if $\rk\Lambda<\infty$ then the decomposition of Proposition
\ref{finite_rank_characterisation200} becomes,
\begin{equation}\label{finite:rank200}
\Lambda=\biggl(\cdots\biggl(\biggl(\widehat{\Lambda}_v
\bigvee_{\Theta_1} T_1\biggr)
\bigvee_{\Theta_2} T_2\biggr)
\cdots\biggr)\bigvee_{\Theta_k} T_k\biggr),
\end{equation}
with $\widehat{\Lambda}_v$ finite, the $\Theta_i$ single vertices, the
$\Theta_i\hookrightarrow\widehat{\Lambda}_v$,
and the images $\Theta_i\hookrightarrow T_i$ incident with a single arc.
Moreover, if $\Lambda\rightarrow\Delta$ is a covering with $\Delta$
single vertexed and $\Lambda$ of finite rank, 
then by Proposition \ref{topology:coverings:result400}(i), each 
tree $T_i$ realizes an embedding $\R\hookrightarrow\Lambda$ 
of the real line in $\Lambda$, and as the
spine is finite, the trees are thus paired
$$
\begin{pspicture}(0,0)(14,2)
\rput(3,0){
\rput(4,1){\BoxedEPSF{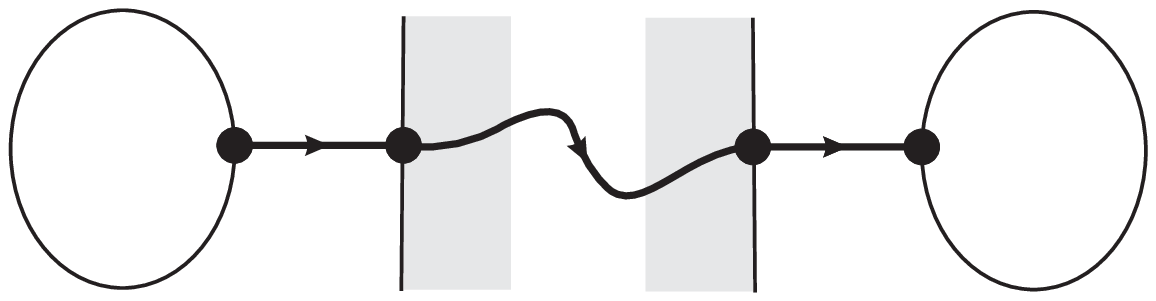 scaled 650}}
\rput(4,1.7){$\widehat{\Lambda}$}
\rput(3.75,.9){$\gamma$}
\rput(2.2,1.3){$e_1$}\rput(5.6,1.3){$e_2$}
\rput(1,1){$T_{e_1}$}\rput(7,1){$T_{e_2}$}
}
\rput(14,1){$(\ddag)$}
\end{pspicture}
$$
with the $e_i$ (and indeed all the edges in the path 
$\R\hookrightarrow\Lambda$) 
in the same fiber of the covering. This pairing will play a key role
in \S \ref{section:pullbacks}.

\begin{corollary}\label{finiterank:result500}
Let $\Lambda\rightarrow\Delta$ be a covering with $\Delta$ 
single vertexed having
non-empty edge set and $\rk\Lambda<\infty$. Then 
$\text{deg}(\Lambda\rightarrow\Delta)<\infty$ 
if and only if $\Lambda=\widehat{\Lambda}_v$.
\end{corollary}

\begin{proof}
If $\Lambda$ is more than $\widehat{\Lambda}_v$ then one of the trees $T_i$ in 
the decomposition (\ref{finite:rank200}) is non trivial and by Proposition
\ref{topology:coverings:result400}(i) we get a real line subgraph
$\R\hookrightarrow\Lambda$, with image in the fiber of an edge, contradicting
the finiteness of the degree. The converse follows from Corollary 
\ref{invariants:result400}.
\qed
\end{proof}

\begin{proposition}
\label{finiterank:result600}
Let $\Lambda\rightarrow\Delta$ be a covering with (i). $\rk\Delta>1$, (ii). 
$\rk\Lambda<\infty$, and (iii). for any intermediate covering 
$\Lambda\rightarrow\Gamma\rightarrow\Delta$ we have $\rk\Gamma<\infty$. Then 
$\text{deg}(\Lambda\rightarrow\Delta)<\infty$.
\end{proposition}

The covering $\R\rightarrow\Delta$ of a single vertexed $\Delta$
of rank $1$ by the real line shows why the $\rk\Delta>1$ condition 
cannot be dropped.

\begin{proof}
By lattice excision we may pass to
the $\Delta$ single vertexed case while preserving (i)-(iii). Establishing
the degree here and passing back to the general $\Delta$ will give the result.
If the degree of the covering
$\Lambda\rightarrow\Delta$ is infinite for $\Delta$ single vertexed, then
by Corollary \ref{finiterank:result500},
in the decomposition (\ref{finite:rank200}) for $\Lambda$, one
of the trees is non-empty and $\Lambda$ has the form
of the graph in Proposition \ref{topology:coverings:result400} with this non-empty
tree the union of the edge $e$ and $\Upsilon_2$. 
Let $\Gamma$ be a graph defined as follows:
take the union of $\Upsilon_1$, the edge $e$ and
$\aa(\R)\cap\Upsilon_2$, where $\aa(\R)$ is the embedding
of the real line given by Proposition \ref{topology:coverings:result400}(i). 
At each vertex of $\aa(\R)\cap\Upsilon_2$
place $\rk\Delta-1$ edge loops:
$$
\begin{pspicture}(0,0)(14,2)
\rput(2,0){
\rput(2,1){\BoxedEPSF{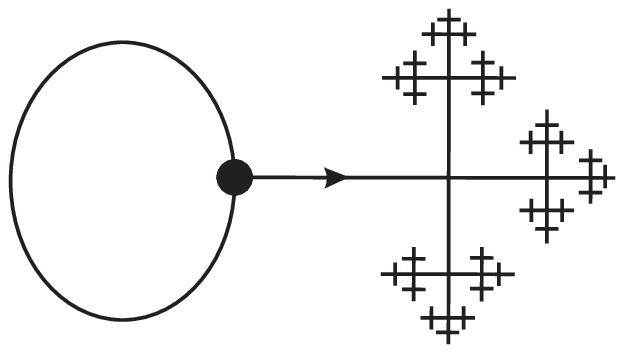 scaled 650}}
\rput(8,1){\BoxedEPSF{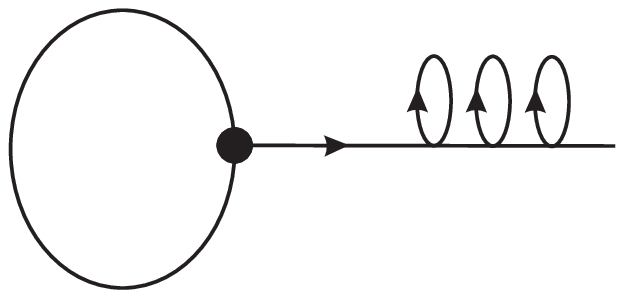 scaled 650}}
\rput(.8,1){$\Upsilon_1$}\rput(6.75,1){$\Upsilon_1$}
\rput(-.2,1.5){$\Lambda$}\rput(5.8,1.5){$\Gamma$}
\rput(9,.7){$\aa(\R)\cap\Upsilon_2$}
\rput(10.2,1.3){$\cdots$}
}
\end{pspicture}
$$
(the picture depicting the $\rk\Delta=2$ case).
Then there is an obvious covering $\Gamma\rightarrow\Delta$ so that
by Proposition \ref{topology:coverings:result400}(ii) we have an intermediate covering
$\Lambda\rightarrow\Gamma\rightarrow\Delta$.
Equally obviously, $\Gamma$ has infinite rank, contradicting (iii). Thus, 
$\text{deg}(\Lambda\rightarrow\Delta)<\infty$.
\qed
\end{proof}

\begin{proposition}
\label{finiterank:result700}
Let $\Psi\rightarrow\Lambda\rightarrow\Delta$
be coverings with $\rk\Lambda<\infty$,
$\Psi\rightarrow\Delta$ Galois, and $\Psi$ not simply connected. Then
$\text{deg}(\Lambda\rightarrow\Delta)<\infty$.
\end{proposition}

The idea is that if the degree is infinite, then
$\Lambda$ has a hanging tree in its 
spine decomposition, and so $\Psi$ does too. But $\Psi$ should look the same
at every point, hence {\em is\/} a tree.

\begin{proof}
Apply lattice excision to 
$\LL(\Psi,\Delta)$, and as $\pi_1(\Psi,u)$ is unaffected by the excision
of trees, we may assume that $\Delta$ is single vertexed.
As $\text{deg}(\Lambda\rightarrow\Delta)$
is infinite, the spine decomposition for $\Lambda$ has an infinite tree,
and $\Lambda$ has the form of Proposition \ref{topology:coverings:result400}. Thus
$\Psi$ does too, by part (iii) of this Proposition, with subgraphs $\Upsilon'_i
\hookrightarrow\Psi$, edge $e'$ and $\Upsilon'_1$ a tree.
Take a closed
reduced path $\gamma$ in $\Upsilon'_2$, and choose a vertex $u_1$ of $\Upsilon'_1$ 
such that the reduced path from $u_1$ to $s(e')$ has at least as many edges
as $\gamma$. Project $\gamma$ via the covering $\Psi\rightarrow
\Delta$ to a closed reduced path, and then lift to $u_1$. The result is
reduced, closed as $\Psi\rightarrow\Delta$ is Galois, and 
entirely contained in the tree $\Upsilon'_1$, hence trivial.
Thus $\gamma$ is also trivial, so that
$\Upsilon'_2$ is a tree and $\Psi$ is simply connected. 
\qed
\end{proof}

\begin{proposition}
\label{finiterank:result800}
Let $\Lambda_u\rightarrow\Delta_v$ be a covering with $\rk\Lambda<\infty$
and $\gamma$ a non-trivial reduced closed path at $v$ lifting to a non-closed path 
at $u$. Then there is an intermediate covering 
$\Lambda_u\rightarrow\Gamma_w\rightarrow\Delta_v$ 
with $\text{deg}(\Gamma\rightarrow\Delta)$
finite and $\gamma$ lifting to a non-closed path at $w$.
\end{proposition}

Stallings shows something very similar \cite[Theorem 6.1]{Stallings83}
starting from a finite immersion rather than a covering. As the proof shows,
the path $\gamma$ in Proposition \ref{finiterank:result800} can be replaced
by finitely many such paths. Moreover, for $T\hookrightarrow\Lambda$ a 
spanning tree, recall that Schreier generators for $\pi_1(\Lambda,u)$
are the homotopy classes of paths through $T$ from $u$ to $s(e)$,
traversing $e$ and traveling back through $T$ to $u_1$,
for $e\in\Lambda\setminus T$. Then
the intermediate $\Gamma$ constructed
has the property that any set of Schreier generators for $\pi_1(\Lambda,u)$
can be extended to a set of Schreier generators for $\pi_1(\Gamma,w)$.

\begin{proof}
If $T\hookrightarrow\Delta$ is a spanning tree and $q:\Delta\rightarrow\Delta/T$
then $\gamma$ cannot be contained in $T$, and so $q(\gamma)$ is non-trivial, closed
and reduced. If the lift of $q(\gamma)$ to $\Lambda/T_i$ is closed then the lift
of $\gamma$ to $\Lambda$ has start and finish vertices that lie in the same
component $T_i$ of $f^{-1}(T)$, mapped isomorphically onto $T$ by the covering, 
and thus implying that
$\gamma$ is not closed. Thus we may apply lattice excision and pass to the
single vertexed case while maintaining $\gamma$ and its properties. Moreover,
the conclusion in this case gives the result in general as closed paths go to closed paths
when excising trees.
If the lift $\gamma_1$ of $\gamma$
at $u$ is not contained in the spine $\widehat{\Lambda}_u$, then its terminal 
vertex lies in a tree $T_{e_i}$ of the spine decomposition $(\ddag)$. By adding
an edge if necessary to $\widehat{\Lambda}_u\cup\gamma_1$, we obtain
a finite subgraph whose coboundary edges are paired, with the edges in each
pair covering the same edge in $\Delta$, as below left:
$$
\begin{pspicture}(0,0)(14,1)
\rput(4.05,.5){$\widehat{\Lambda}_u\cup\gamma_1$}
\rput(10,.35){$\Gamma$}
\rput(4,.5){\BoxedEPSF{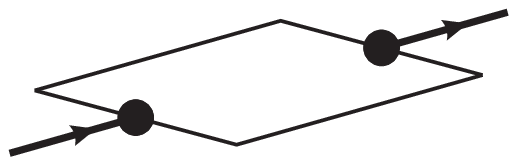 scaled 740}}
\rput(10,.65){\BoxedEPSF{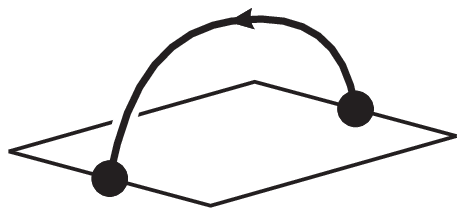 scaled 740}}
\end{pspicture}
$$
(if the lift is contained in the spine,
take $\widehat{\Lambda}_u$ itself).
In any case, let $\Gamma$ be $\widehat{\Lambda}_u\cup\gamma_1$ together with a single 
edge replacing each pair as above right. Restricting the covering 
$\Lambda\rightarrow\Delta$ to $\widehat{\Lambda}_u\cup\gamma_1$ and 
mapping the new edges to the
common image of the old edge pairs gives a finite covering $\Gamma\rightarrow\Delta$,
and hence an intermediate covering
$\Lambda\stackrel{q}{\rightarrow}\Gamma{\rightarrow}\Delta$, with 
$q(\gamma_1)$ non-closed at $q(u)$.
\qed
\end{proof}

For the rest of this section we investigate the rank implications 
of the decomposition (\ref{finite:rank200})
and the pairing $(\ddag)$ in a special case. 
Suppose $\Lambda\rightarrow\Delta$ is a covering 
with $\Delta$ single vertexed, $\rk\Delta=2$, $\Lambda$ non-simply 
connected, and $\rk\Lambda<\infty$. Let $x_i^{\pm 1}, (1\leq i\leq 2)$
be the edge loops of $\Delta$
and fix a spine so we have the decomposition (\ref{finite:rank200}).

An {\em extended spine\/} for such a $\Lambda$ 
is a connected subgraph $\Gamma\hookrightarrow\Lambda$ 
obtained by adding finitely many edges to a spine, so that every
vertex of $\Gamma$ is incident with either zero or three edges
in its coboundary $\delta\Gamma$. It is always possible to find
an extended spine: take the 
union of the spine $\widehat{\Lambda}_u$ and each 
edge $e\in\delta\widehat{\Lambda}_u$ in its coboundary. 
Observe that $\Gamma$ is finite and the decomposition
(\ref{finite:rank200}) gives $\rk\Gamma=\rk\widehat{\Lambda}_u=\rk\Lambda$.
Call a vertex of the extended spine $\Gamma$
{\em interior\/} (respectively {\em boundary\/}) when it is incident
with zero (resp. three) edges in $\delta\Gamma$.

We have the pairing of trees $(\ddag)$ for an extended spine, so that
each boundary vertex $v_1$ is paired with another $v_2$,
$$
\begin{pspicture}(0,0)(14,2)
\rput(3,0){
\rput(4,1){\BoxedEPSF{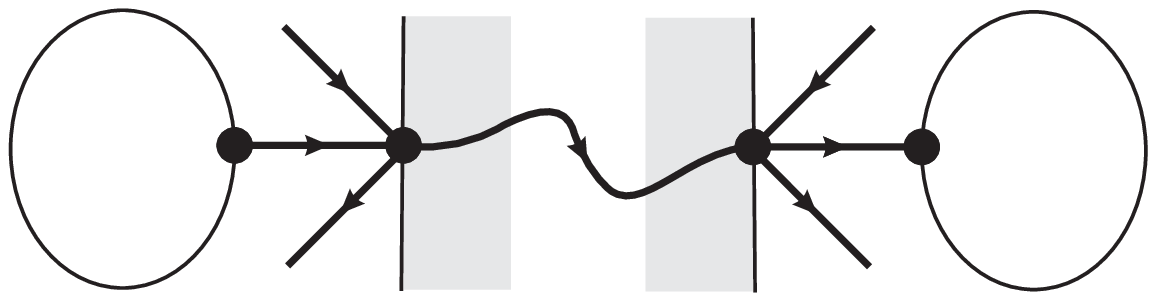 scaled 650}}
\rput(4,1.7){$\Gamma$}
\rput(3.8,.9){$\gamma$}
\rput(2.2,1.25){$e_1$}\rput(5.8,1.25){$e_2$}
\rput(3.1,1.25){$v_1$}\rput(4.9,1.25){$v_2$}
\rput(1,1){$T_{e_1}$}\rput(7,1){$T_{e_2}$}
}
\rput(14,1){$(*)$}
\end{pspicture}
$$
with $e_1,e_2$ and all the edges in the path $\gamma=\aa(\R)\cap\Gamma$ 
covering an edge loop $x_i\in\Delta$. Call this an {\em $x_i$-pair\/},
($i=1,2$).

For two $x_i$-pairs (fixed $i$), the
respective $\gamma$ paths share no vertices in common, for otherwise
there would be two distinct edges covering the same $x_i\in\Delta$ starting
at such a common vertex. 
Moreover, $\gamma$ must contain vertices of $\Gamma$ apart from
the two boundary vertices $v_1,v_2$, otherwise $\Lambda$ 
would be simply connected.
These other vertices are incident with at least two edges of $\gamma\in\Gamma$, 
hence at most $2$ edges of the coboundary $\delta\Gamma$, and thus must be 
interior.

\begin{lemma}\label{finiterank:result1000}
If $n_i, (i=1,2)$, is the number of $x_i$-pairs in an extended spine $\Gamma$ for
$\Lambda$, 
then the number of interior vertices is at least $\sum n_i$.
\end{lemma}

\begin{proof}
The number of interior vertices is $|V_\Gamma|-2\sum n_i$ and the number of
edges of $\Gamma$ is $4(|V_\Gamma|-2\sum n_i)+2\sum n_i$, hence 
$\rk\Gamma-1=|V_\Gamma|-3\sum n_i$. 
As $\Lambda$ is not simply connected,
$\rk\Lambda-1=\rk\Gamma-1\geq 0$, thus $|V_\Gamma|-2\sum n_i\geq
\sum n_i$ as required.
\qed
\end{proof}

The lemma is not true in the 
case $\rk\Delta>2$.
It will be helpful in \S \ref{section:pullbacks} 
to have a pictorial description of the quantity
$\rk-1$ for our graphs. To this end, a {\em checker\/} is a small plastic
disk, as used in the eponymous boardgame (called {\em draughts\/}
in British English). We place black checkers on some of the
vertices of an extended spine $\Gamma$ according to the 
following scheme: place black checkers on all the interior vertices
of $\Gamma$; for each $x_1$-pair in (*), take the interior vertex
on the path $\gamma$ that is closest to $v_1$ (ie: is the terminal vertex of the
edge of $\gamma$ whose start vertex is $v_1$) and {\em remove\/} its
checker; for each
$x_2$-pair, we can find, by Lemma \ref{finiterank:result1000}, an interior
vertex with a checker still on it. Choose such a vertex and remove its checker
also.
We saw in the proof of Lemma \ref{finiterank:result1000} that $\rk\Lambda-1=
\rk\Gamma-1$ is equal to the number of interior vertices of $\Gamma$, less
the number of $x_i$-pairs $(i=1,2)$. Thus,

\begin{lemma}\label{finiterank:whitevertices}
With black checkers placed on the vertices of an extended spine for $\Lambda$ 
as above, the number of black checkers is $\rk\Lambda-1$.
\end{lemma}

From now on we will only use the extended spine obtained by adding the 
coboundary edges to some fixed spine $\widehat{\Lambda}_u$.

Let $p:\Lambda_u\rightarrow\Delta_v$ be a covering with $\rk\Delta=2$,
$\rk\Lambda<\infty$ and $\Lambda$ not simply connected. A spanning tree
$T\hookrightarrow\Delta$ induces a covering $\Lambda/T_i\rightarrow\Delta/T$
with $\Delta/T$ single vertexed. Let $\HH(\Lambda_u\rightarrow\Delta_v)$ be the
number of vertices of the spine of $\Lambda/T_i$ at $q(u)$ and 
$n_i(\Lambda_u\rightarrow\Delta_v)$ the number of $x_i$-pairs in the 
extended spine. The isomorphism class of $\Lambda/T_i$ and the spine
are independent of the spanning tree $T$, hence the quantities
$\HH(\Lambda_u\rightarrow\Delta_v)$ and $n_i(\Lambda_u\rightarrow\Delta_v)$
are too.

\section{Pullbacks}\label{section:pullbacks}

Let $p_i:\Lambda_i:=\Lambda_{u_i}\rightarrow\Delta_v, (i=1,2)$ be coverings and 
$\Lambda_{1}\prod_\Delta\Lambda_{2}$ their (unpointed) pullback. 
If $\widehat{\Lambda}_{u_i}$ is the spine at $u_i$ then we can restrict
the coverings to maps $p_i:\widehat{\Lambda}_{u_i}\rightarrow\Delta_v$
and form the pullback $\widehat{\Lambda}_{u_1}\prod_\Delta\widehat{\Lambda}_{u_2}$.

\begin{proposition}[spine decomposition of pullbacks]
\label{pullbacks:spinedecomposition}
The pullback $\Lambda=\Lambda_{1}\prod_\Delta\Lambda_{2}$ 
has a wedge sum
decomposition 
$\Lambda=(\widehat{\Lambda}_{u_1}\prod_\Delta\widehat{\Lambda}_{u_2})
\bigvee_\Theta\Phi$,
with $\Phi$ a forest 
and no two vertices of the image of 
$\Theta\hookrightarrow\Phi$ lying in the same component.
\end{proposition}

\begin{proof}
Let $\Lambda_i=\widehat{\Lambda}_{u_i}\bigvee_{\Theta_i}\Phi_i, (i=1,2)$ be the 
spine 
decomposition, $t_i:\Lambda_1\prod_\Delta\Lambda_2\rightarrow\Lambda_i, (i=1,2)$
the coverings provided by the pullback and
$\Omega$ a connected component of the pullback. If 
$\Omega\cap(\widehat{\Lambda}_{u_1}\prod_\Delta\widehat{\Lambda}_{u_2})
=\varnothing$, then a reduced closed path $\gamma\in\Omega$ must map via one of 
the $t_i$ to a closed path in the forest $\Phi_i$. As the images under 
coverings of
reduced paths are reduced, $t_i(\gamma)$ must contain a spur which can be lifted
to a spur in $\gamma$. Thus $\Omega$ is a tree. 

Otherwise choose 
a vertex $w_1\times w_2$ in 
$\Omega\cap(\widehat{\Lambda}_{u_1}\prod_\Delta\widehat{\Lambda}_{u_2})$ and
let $\Gamma$ be the connected component of this intersection containing
$w_1\times w_2$. If $\gamma$ a reduced closed path at $w_1\times w_2$ then
$t_i(\gamma), (i=1,2)$ is a reduced closed path at 
$w_i\in\widehat{\Lambda}_{u_i}$, hence 
$t_i(\gamma)\in\widehat{\Lambda}_{u_i}$ and thus 
$\gamma\in\widehat{\Lambda}_{u_1}\prod_\Delta\widehat{\Lambda}_{u_2}$.
Applying Proposition \ref{finite_rank_characterisation200}, we have $\Omega$
a wedge sum of $\Gamma$ and a forest of the required form.
\qed
\end{proof}


\begin{corollary}[Howsen-Stallings]
\label{pullbacks:result200}
Let $p_i:\Lambda_i\rightarrow\Delta, (i=1,2),$ be coverings with 
$\rk\Lambda_i<\infty$ and $u_1\times u_2$ a vertex of their pullback. Then 
$\rk(\Lambda_1\prod_\Delta\Lambda_2)_{u_1\times u_2}<\infty.$
\end{corollary}

\begin{proof}
The component $\Omega$ of the pullback containing $u_1\times u_2$ is either a tree
or the wedge sum of a finite graph and a forest as described in Proposition
\ref{pullbacks:spinedecomposition}. Either case gives the result.
\qed
\end{proof}

The remainder of this section is devoted to a proof of an estimate for the
rank of the pullback of finite rank graphs in a special case. Let 
$p_j:\Lambda_j:=\Lambda_{u_j}\rightarrow\Delta_v, (j=1,2)$ be coverings
with $\rk\Delta=2$, $\rk\Lambda_j<\infty$ and the
$\Lambda_j$ not simply connected. 
Let $\HH_j:=\HH(\Lambda_{u_j}\rightarrow\Delta_v)$ and 
$n_{ji}:=n_i(\Lambda_{u_j}\rightarrow\Delta_v)$ be as at the end of 
\S \ref{section:invariants}.

\begin{theorem}\label{pullback:rankestimate}
For $i=1,2$,
$$
\sum_\Omega (\rk\Omega-1)\leq \prod_j(\rk\Lambda_j-1)
+\HH_1\HH_2-(\HH_1-n_{1i})(\HH_2-n_{2i}),
$$
the sum over all non simply connected components $\Omega$ of the pullback
$\Lambda_1\prod_\Delta\Lambda_2$.
\end{theorem}

\begin{proof}
Lattice excision and the definition of the $\HH_j$ and $n_{ji}$ allow us to pass
to the $\Delta$ single vertexed case. 
Suppose then that $\Delta$ has edge loops 
$x_i^{\pm 1}, (1\leq i\leq 2)$ at the vertex $v$, and
extended spines 
$\widehat{\Lambda}_{u_j}\hookrightarrow\Gamma_j\hookrightarrow\Lambda_j$.
The covering $p_j:\Lambda_j\rightarrow\Delta_v$ can be restricted to maps
$\Gamma_j\rightarrow\Delta_v$ and $\widehat{\Lambda}_{u_j}\rightarrow\Delta_v$,
and we form the three resulting pullbacks 
$\Lambda_1\prod_\Delta\Lambda_2$, $\Gamma_1\prod_\Delta\Gamma_2$
and $\widehat{\Lambda}_{u_1}\prod_\Delta\widehat{\Lambda}_{u_2}$, with
$$
\widehat{\Lambda}_{u_1}\prod_\Delta\widehat{\Lambda}_{u_2}\hookrightarrow
\Gamma_1\prod_\Delta\Gamma_2\hookrightarrow
\Lambda_1\prod_\Delta\Lambda_2,
$$ 
and $t_j:\Lambda_1\prod_\Delta\Lambda_2
\rightarrow\Lambda_j$ the resulting covering maps. 

Place black checkers on the vertices of the extended spines $\Gamma_j$ as in 
\S \ref{section:invariants} and place a black checker on
a vertex $v_1\times v_2$ of $\Gamma_1\prod_\Delta\Gamma_2$ 
precisely when both $t_j(v_1\times v_2)\in\Gamma_j, (j=1,2)$ have black checkers
on them.
By Lemma \ref{finiterank:whitevertices}, and the construction of the
pullback for $\Delta$ single vertexed, we get the number of vertices
in $\Gamma_1\prod_\Delta\Gamma_2$ with black checkers is equal to 
$\prod(\rk\Lambda_j-1)$.

Let $\Omega$ be a non simply connected component of the pullback 
$\Lambda_1\prod_\Delta\Lambda_2$ and 
$\Upsilon=\Omega\cap(\Gamma_1\prod_\Delta\Gamma_2)$.
If $v_1\times v_2$ is the start vertex of
at least one edge in the coboundary $\delta\Upsilon$,
then at least one of the $v_j$ must be incident with at least one,
hence three, edges of the coboundary $\delta\Gamma_j$. 
Lifting these three via
the covering $t_j$ to $v_1\times v_2$ gives at least three edges 
starting at $v_1\times v_2$ in the coboundary
$\delta\Upsilon$. Four coboundary edges starting here
would mean that $\Omega$ was simply connected,
hence every vertex of $\Upsilon$
is incident with either
zero or three coboundary edges. 

We can thus extend the interior/boundary terminology of \S \ref{section:invariants}
to the vertices of 
$\Upsilon$, and observe that a vertex of $\Upsilon$ covering,
via either of the $t_j$, a boundary vertex $\in\Gamma_j$, must itself be a
boundary vertex.
The upshot is that $\Upsilon$ is an extended spine 
in $\Omega$ and by Proposition \ref{pullbacks:spinedecomposition}, $\rk\Omega-1=
\rk\Upsilon-1$. Now place {\em red\/} checkers on the vertices of $\Upsilon$  
as in \S \ref{section:invariants} and do this for each non-simply connected 
component $\Omega$. The number of red checkered vertices is 
$\sum_\Omega (\rk\Omega-1)$.

The result is that $\Gamma_1\prod_\Delta\Gamma_2$ has vertices 
with black checkers, vertices with red checkers, vertices with red checkers sitting
on top of black checkers, and vertices that are completely uncheckered. Thus,
$$
\sum_\Omega (\rk\Omega-1)\leq \prod(\rk\Lambda_j-1)+N,
$$
where $N$ is the number of vertices of $\Gamma_1\prod_\Delta\Gamma_2$ that 
have a red checker but no black checker.

\parshape=11 0pt\hsize 0pt\hsize 0pt\hsize 0pt\hsize 0pt\hsize 
0pt\hsize 0pt\hsize 0pt\hsize 
0pt\hsize 0pt.6\hsize 0pt.6\hsize 
It remains then to estimate the number of these ``isolated'' red checkers. Observe
that a vertex of $\Gamma_1\prod_\Delta\Gamma_2$ has no black checker precisely
when it lies in the fiber, via at least one of the $t_j$, of a checkerless vertex
in $\Gamma_j$. Turning it around, we investigate the fibers of the checkerless 
vertices of both $\Gamma_j$. 
Indeed, in an $x_1$-pair, 
the vertices $v_1,v_2$ and $u$ are checkerless, while $v_1,v_2$ are also checkerless
in an $x_2$-pair. We claim that no vertex in the fiber, via $t_j$, of these
five has a red checker. A vertex of $\Upsilon$ in the fiber of the boundary
vertices $v_1,v_2$ is itself a boundary vertex, hence contains no red checker.
If $v_1\times v_2\in\Upsilon$ is in the fiber of $u$ and is a boundary vertex of 
$\Upsilon$ then it carries no red checker either. If instead 
$v_1\times v_2$ is an interior vertex
then the lift to $v_1\times v_2$ of $e^{-1}$ cannot be in the coboundary
$\delta\Upsilon$, hence the terminal vertex of this lift is in $\Upsilon$ also
and covers $v_1$. Thus, this terminal vertex is a boundary vertex for an $x_1$-pair
of $\Upsilon$, and $v_1\times v_2$ is the interior vertex from which a red
checker is removed for this pair.
\vadjust{\hfill\smash{\lower 8mm
\llap{
\begin{pspicture}(0,0)(5,2)
\rput(-1.75,-.1){
\rput(4,1){\BoxedEPSF{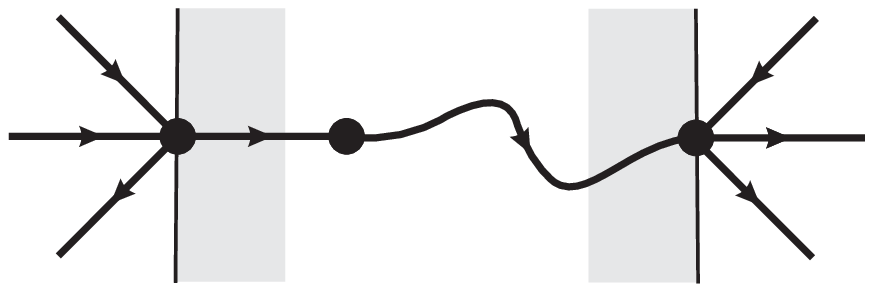 scaled 650}}
\rput(4,1.7){$\Gamma_j$}
\rput(2.8,.85){$e$}
\rput(2.55,1.25){$v_1$}\rput(5.4,1.25){$v_2$}
\rput(3.6,1.25){$u$}
}
\end{pspicture}
}}}\ignorespaces

\parshape=3 
0pt.6\hsize 0pt.6\hsize 0pt\hsize 
The only remaining checkerless vertices of the $\Gamma_j$ unaccounted for are
those interior vertices chosen for each $x_2$-pair.
Let $S_1=\{u_1,\ldots,u_{n_{12}}\}\subset\Gamma_1$ and 
$S_2=\{w_1,\ldots,w_{n_{22}}\}\subset\Gamma_2$
be these sets of vertices. The result of the discussion above is that if
$v_1\times v_2$ has an isolated red checker then it must be contained in
$(S_1\times V_{\Gamma_2})\cup(V_{\Gamma_1}\times S_2)$, 
the vertices of $\Gamma_1\prod_\Delta\Gamma_2$ in the fiber of
a $u_i$ or a $w_i$. If $u_i\times y\in S_1\times V_{\Gamma_2}$ with $y$ 
a boundary vertex of $\Gamma_2$, then $u_i\times y$ is a boundary vertex of
$\Gamma_1\prod_\Delta\Gamma_2$, hence has no red checker. 
Similarly a $x\times w_i$ with $x$ a boundary vertex of $\Gamma_1$
has no red checker, and so $N$ is at most
the number of vertices in the set $(S_1\times V_2)\cup(V_1\times S_2)$, with
$V_i$ the vertices of the spine $\widehat{\Lambda}_{u_i}$. As $S_i\subset V_i$, the
two sets in this union intersect in $S_1\times S_2$, so we have
$$
N\leq |S_1\times V_2|+|V_1\times S_2|-|S_1\cap S_2|=
n_{12}\HH_2+n_{22}\HH_1-n_{12}n_{22},
$$
hence the result for $i=2$.
Interchanging the checkering scheme for the $x_i$-pairs
gives the result for $i=1$.
\qed
\end{proof}

\section{Free groups and the topological dictionary}\label{free}

A group $F$ is {\em free of rank $\rk F$\/} if and only if it is 
isomorphic to the fundamental group of a connected graph of rank $\rk F$.
If\/ $\Gamma_1,\Gamma_2$ are connected graphs with 
$\pi_1(\Gamma_1,v_1)\cong\pi_1(\Gamma_2,v_2)$, then
$H_1(\Gamma_1)\cong H_1(\Gamma_2)$ and thus 
$\rk\Gamma_1=\rk\Gamma_2$.

The free groups so defined are of course the standard free groups and the 
rank is the usual rank
of a free group. At this stage we appeal to the existing (algebraic) theory 
of free groups, and in particular, 
that by applying Nielsen transformations, a set of generators for a free 
group 
can be transformed into a set of 
free generators whose cardinality is no greater. Thus, a finitely generated 
free group has finite rank (the converse
being obvious). From now on we use the (topologically more tractable) notion
of finite rank as a synonym for finitely generated.

Let $F$ be a free group and
$\varphi:F\rightarrow\pi_1(\Delta,v)$ an isomorphism for 
$\Delta$ connected. We call $\vphi$ a topological realization, and
the ``topological dictionary'' is the loose term used to describe the
correspondence between algebraic properties of $F$ and topological
properties of $\Delta$.
The non-abelian $F$ correspond to the
$\Delta$ with $\rk\Delta>1$.
A subgroup $A\subset F$ corresponds to a 
covering $f:\Lambda_u\rightarrow\Delta_v$
with $f_*\pi_1(\Lambda,u)=\varphi(A)$, and hence $\rk A=\rk\Lambda$
($f_*$ is the homomorphism induced by $p$ using the functorality of $\pi_1$).
Thus
finitely generated subgroups correspond to finite rank $\Lambda$ and normal subgroups
to Galois coverings. Inclusion relations between subgroups correspond to 
covering relations, indices of subgroups to degrees of coverings, trivial
subgroups to simply connected coverings, conjugation to change of basepoint, and
so on.

Applying the topological dictionary to the italicised results below we recover some
classical facts (see also \cite{Servatius83,Stallings83}).

\begin{enumerate}
\item \cite{Greenberg60,Karrass69}: If a finitely generated subgroup $A$ of a 
non-abelian free group $F$ 
is contained in no subgroup of infinite rank, then $A$ has finite index in $F$;
{\em Proposition \ref{finiterank:result600}}.
\item \cite{Greenberg60}: 
If a finitely generated subgroup $A$ of a free group 
$F$ contains a non-trivial normal subgroup of $F$, then it has finite index in $F$;
{\em Proposition \ref{finiterank:result700}}.
\item \cite{Burns69,Hall49}: Let $F$ be a free group, $X$ a finite subset of $F$, 
and $A$ a finitely
generated subgroup of $F$ disjoint from $X$. Then $A$ is a free
factor of a group $G$, of finite index in $F$ and disjoint from $X$;
{\em Proposition \ref{finiterank:result800}} (and the comments following it).
\item \cite{Howsen54}: If $A_1,A_2$ are finitely generated subgroups 
of a free group $F$, then the intersection of conjugates
$A_1^{g_1}\cap A_2^{g_2}$ is finitely generated for any
$g_1,g_2\in F$;
{\em Corollary \ref{pullbacks:result200}}.
\end{enumerate}

If $\Delta$ is a graph, $\rk\Delta=2$, and $A\subset F=\pi_1(\Delta,v)$, then we define
$\HH(F,A):=\HH(\Lambda_u\rightarrow\Delta_v)$ and 
$n_i(F,A):= n_i(\Lambda_u\rightarrow\Delta_v)$, where
$f:\Lambda_u\rightarrow\Delta_v$ is the covering with $f_*\pi_1(\Lambda,u)=A$. 
For an arbitrary free group $F$ realized via
$\varphi:F\rightarrow\pi_1(\Delta,v)$,
define $\HH^\varphi(F,A)$ and $n^\varphi_i(F,A)$ to be
$\HH(\varphi(F),\varphi(A))$ and $n_i(\varphi(F),\varphi(A))$.

The appearance of $\varphi$ in the notation is meant to indicate that these 
quantities, unlike rank, are realization dependent. This can be both a strength
and a weakness. A weakness because it seems desirable for algebraic statements to
involve only algebraic invariants, and a strength if we have the freedom to choose
the realization, especially if more interesting results are obtained when this
realization is not the ``obvious'' one.

For example, if $F$ is a free group with free generators $x$ and $y$, and $\Delta$ is
single vertexed with two edge loops whose homotopy classes are $a$ and $b$, then
the subgroup $A=\langle xy\rangle\subset F$ corresponds to the $\Lambda$ below
left under the obvious realization $\varphi_1(x)=a,\varphi_1(y)=b$, and
to the righthand graph via $\varphi_2(x)=a,\varphi_2(y)=a^{-1}b$:
$$
\begin{pspicture}(0,0)(12,3)
\rput(9.5,1.5){\BoxedEPSF{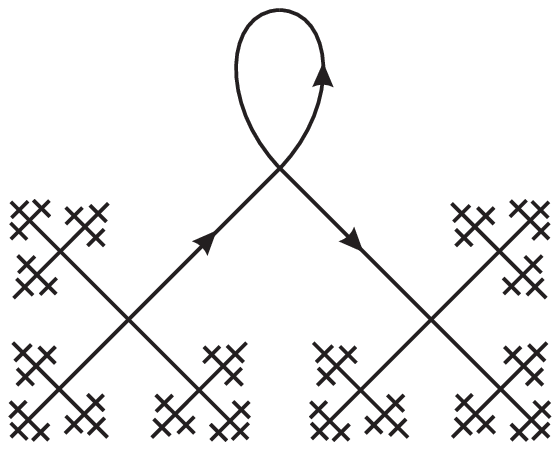 scaled 500}}
\rput(3.5,1.5){\BoxedEPSF{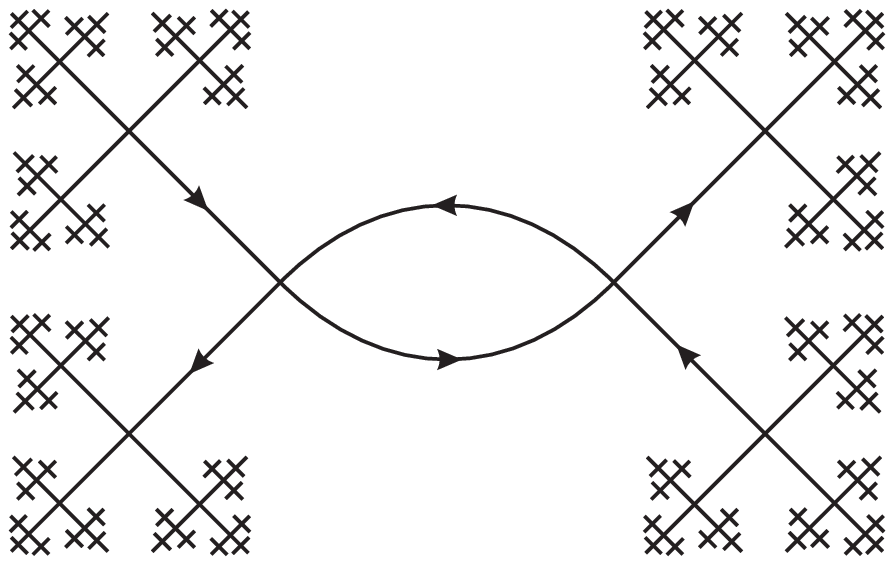 scaled 500}}
\end{pspicture}
$$
Thus, $\HH^{\varphi_1}(F,A)=2,n^{\varphi_1}_{i}(F,A)=1, (i=1,2)$, 
whereas $\HH^{\varphi_2}(F,A)=1,n^{\varphi_2}_{1}(F,A)=1,n^{\varphi_2}_{2}(F,A)=0$.

We now apply the topological dictionary to Theorem \ref{pullback:rankestimate}.
Let $\varphi:F\rightarrow\pi_1(\Delta,v)$, $A_j\subset F, (j=1,2)$, finitely generated
non-trivial subgroups, and 
$f_j:\Lambda_{u_j}\rightarrow\Delta_v, (j=1,2)$ coverings with
$\varphi(A_j)={f_{j}}_*\pi_1(\Lambda,u_j)$.
Each non simply-connected component $\Omega$ of the pullback corresponds to some
non-trivial intersection of conjugates $A_1^{g_1}\cap A_2^{g_2}$. 
As observed in \cite{Neumann90}, these in turn
correspond to the conjugates $A_1\cap A_2^g$ for $g$ from a set of
double coset representatives for $A_2\backslash F/ A_1$.

\begin{theorem}
\label{algebraic:shn}
Let $F$ be a free group of rank two and $A_j\subset F, (j=1,2)$, finitely generated
non-trivial subgroups. Then for any realization $\varphi:F\rightarrow\pi_1(\Delta,v)$
and $i=1,2$, 
$$
\sum_g (\rk(A_1\cap A_2^g)-1)\leq 
\prod_j(\rk A_j-1)
+\HH_1\HH_2-(\HH_1-n_{1i})(\HH_2-n_{2i}),
$$
the sum over all double coset representatives $g$ for $A_2\backslash F/ A_1$
with $A_1\cap A_2^g$ non-trivial, and where $\HH_j=\HH^\varphi(F,A_j)$ and
$n_{ji}=n^\varphi_i(F,A_j)$.
\end{theorem}

This theorem should be viewed in the context of attempts 
to prove the so-called {\em strengthened Hanna Neumann conjecture\/}: namely,
if $A_j, (j=1,2)$ are finitely
generated, non-trivial, subgroups of an arbitrary free group $F$, then 
$$
\sum_g (\rk(A_1\cap A_2^g)-1)\leq 
\prod_j(\rk A_j-1)+\varepsilon,
$$
the sum over all double coset representatives $g$ for $A_2\backslash F/ A_1$
with $A_1\cap A_2^g$ non-trivial,
where the conjecture is that $\varepsilon$ is zero, while in the existing
results, it is an error term having a long history. 
A selection of estimates for $\ve$, in chronological order is,
$(\rk A_1-1)(\rk A_2-1)$ \cite{Neumann56},
$\max\{(\rk A_1-2)(\rk A_2-1),(\rk A_1-1)(\rk A_2-2)\}$, \cite{Burns69},
$\max\{(\rk A_1-2)(\rk A_2-2)-1,0\}$, \cite{Tardos96} and
$\max\{(\rk A_1-3)(\rk A_2-3),0\}$ \cite{Dicks01}
(the original, unstrengthened conjecture \cite{Neumann56} involved just the
intersection of the two subgroups, rather than their conjugates, and
the first two expressions for $\varepsilon$ were proved in this restricted
sense; the strengthened version was formulated in \cite{Neumann90}, and the 
H. Neumann and Burns estimates for $\varepsilon$ were improved to the 
strengthened case there).
Observe that as the join
$\langle A_1,A_2\rangle$ of two finitely generated subgroups 
is finitely generated, and every finitely generated free
group can be embedded as a subgroup of the free group of rank two, we may replace
the ambient free group in the conjecture with the free group of rank two.

It is hard to make a precise comparison between the $\varepsilon$ provided by 
Theorem \ref{algebraic:shn}
and those above. Observe that if $A_j\subset F$, with
$F$ free of rank two, then with respect to a topological realization we have
$\rk A_j=\HH_j-(n_{j1}+n_{j2})+1$. It is straightforward to find infinite families
$A_{1k},A_{2k}\subset\pi_1(\Delta,v), (k\in\Z^{>0})$, for which the error term in
Theorem \ref{algebraic:shn} is less than those listed above for all but finitely
many $k$, or even for which the strengthened Hanna Neumann conjecture is true
by Theorem \ref{algebraic:shn}, for instance,
$$
\begin{pspicture}(0,0)(14.5,2)
\rput(6,1){\BoxedEPSF{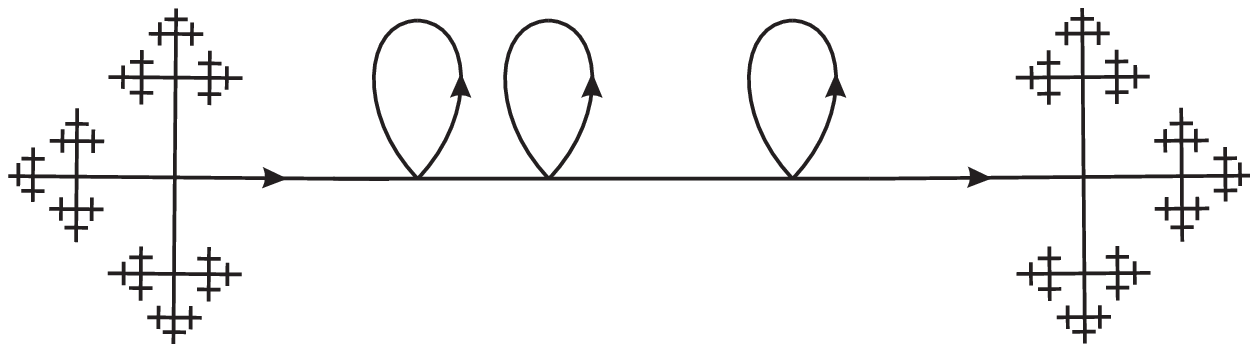 scaled 500}}
\rput{270}(5.9,.8){$\left.\begin{array}{c}
\vrule width 0 mm height 22 mm depth 0 pt\end{array}\right\}$}
\rput(5.9,.4){$k$ edge loops}
\rput(1.8,1.5){$A_{1k}=A_{2k}=$}\rput(6.2,1.5){$\ldots$}
\rput(1,-.2){\rput(10.5,1.75){$\HH_i=k$}
\rput(10.5,1){$n_{11}=n_{21}=0$}
\rput(10.5,.5){$n_{12}=n_{22}=1$}}
\end{pspicture}
$$
but where the error terms above are quadratic in $k$.


{}

\end{document}